\documentclass{amsart}
\usepackage{amsmath, amssymb}
\usepackage{amsfonts}
\usepackage[arrow,matrix,curve,cmtip,ps]{xy}

\usepackage{amsthm}

\allowdisplaybreaks

\newtheorem{theorem}{Theorem}[section]
\newtheorem{lemma}[theorem]{Lemma}

\newtheorem*{theorem*}{Theorem}
\theoremstyle{remark}
\newtheorem{remark}[theorem]{Remark}
\newtheorem{definition}[theorem]{Definition}

\newcommand{\thmref}[1]{Theorem~\ref{#1}}

\newcommand{\lemref}[1]{Lemma~\ref{#1}}

\numberwithin{equation}{section}

\newcommand{\C}{\mathbb{C}}

\newcommand{\B}{\mathcal{B}}

\newcommand{\V}{\mathcal{V}}

\newcommand{\go}{G^0}
\newcommand{\cc}{C_{c}}

\newcommand{\cs}{\ensuremath{C^{*}}}
\newcommand{\supp}{\operatorname{supp}}
\newcommand{\csg}[1]{\cs(G|_{#1})}

\newcommand{\co}{C_{0}}

\renewcommand{\H}{\mathcal{H}}

\newcommand{\Ind}{\operatorname{Ind}}

\begin{document}


\title{CCR and GCR groupoid $\cs$-algebras}

\author{Lisa Orloff Clark 
}

\curraddr{Dept of Mathematical Sciences\\ Susquehanna University\\
Selinsgrove\\ PA 17870\\ USA} \email{clarklisa@susqu.edu}

\date{\today}
\subjclass{46L05,46L35}
\keywords{Locally Compact groupoid, $\cs$-algebra, Hilbert Module}
\begin{abstract}

Suppose $G$ is a second countable, locally compact,
Hausdorff groupoid with a fixed left Haar system.
Let $\go/G$ denote the orbit space of $G$ and $\cs(G)$ denote
the groupoid $\cs$-algebra.  Suppose that the isotropy groups of
$G$ are amenable.  We show that $\cs(G)$ is CCR if and only if $\go/G$ is a $T_1$ topological space and all of the isotropy groups are CCR.  We also show that
$\cs(G)$ is GCR if and only if $\go/G$ is a $T_0$ topological space and all of the isotropy groups are GCR.
\end{abstract}

\maketitle

\section{Intorduction}
In this paper we generalize two major classification theorems for transformation group $\cs$-algebras  to theorems that classify groupoid $\cs$-algebras.  We will denote a transformation group by $\cs(H,X)$ where $H$ is a locally compact group acting continuously on the left of the locally compact space $X$.
Specifically, in \cite{gootman}, Elliot Gootman showed  the following:
\begin{theorem}
\label{tgthm4}
Suppose $H$ and $X$ are both second countable.
Then $\cs(H,X)$ is GCR  if and only if the orbit space is $T_0$ and every stability
group is GCR.
\end{theorem}

Dana Williams used a different approach and, in addition to the GCR case, considered the case for CCR transformation
group $\cs$-algebras.  In \cite{danasthesis} he proved the theorem below.

\begin{theorem}
\label{tgthm5}
Suppose that $H$ and $X$ are both second countable.
Suppose also that at every point of discontinuity $y$
of the map $x \mapsto S_x$, the stability group $S_y$ is amenable,
then $\cs(H,X)$ is CCR if and only if the orbit space is $T_1$ and the stability
groups are CCR.
\end{theorem}

We will prove analogous results to \thmref{tgthm4} and
\thmref{tgthm5} for groupoid $\cs$-algebras.   Like Williams'
results, we must still assume that the isotropy groups are amenable.
Suppose $G$ is a locally compact, Hausdorff, second countable
groupoid with a fixed left Haar system and $\cs(G)$ is the
associated groupoid $\cs$-algebra.  We show that:

\begin{theorem}
The groupoid $\cs$-algebra
$\cs(G)$ is CCR if and only if the orbit space is $T_1$ and  the isotropy groups are CCR.
\end{theorem}

\begin{theorem}
The groupoid $\cs$-algebra
$\cs(G)$ is GCR if and only if the orbit space is $T_0$ and the isotropy groups are GCR.
\end{theorem}

\thmref{tgthm4} and \thmref{tgthm5} have been generalized to
principal groupoids in \cite{mypaper}. The key to the proof of the
CCR result in the principal case has two major components.  The
first is that every irreducible representation of $\cs(G)$ factors
through $\cc(G_{\overline{[u]}})$ for some $u \in \go$
\cite[Corollary~3.5]{mypaper}.  This result does not require $G$ to
be principal so we can use it in the more general setting. The
second component used in the proof of the CCR result is a continuous
injection from the orbit space of the groupoid $\go/G$ to the
spectrum of the associated groupoid $\cs$-algebra $\cs(G)^\wedge$.
In \cite{mypaper}, this continuous injection is defined explicitly.
However, the map there requires $G$ to be principal in order to be a
well-defined function. Our goal in this paper is to show that there
still exists a continuous injection from the orbit space to the
spectrum in the more general setting. Once we do this, the CCR
result will be easy to prove. The GCR result will also follow
quickly because, as in the principal case, we can invoke
\cite[Proposition~5.1]{mypaper}.

 The authors of \cite{ctstraceIII} have also defined a map from $\go/G$ to
$\cs(G)^\wedge$.  Their definition assumes the groupoid has abelian
isotropy. What we will do is show that we can modify this definition
and show that this map is a well-defined injection without the
assumption that the groupoid has abelian isotropy.  This is highly
non-trivial as it involves possibly discontinuous isotropy. In fact,
we must invoke Renault's Disintegration theorem
\cite[Theorem~4.2]{renault2} to show our map gives us a bounded
operator, which is not needed in the abelian case.  Unfortunately,
to show this map is continuous, we must still require the isotropy
groups be amenable.  In \thmref{tgthm4}, Gootman does not require
the amenability condition. Thus, we expect it to be superfluous in
the groupoid case, but, with these techniques, we have been
unsuccessful in eliminating it.

To define a map from $\go/G$ to $\cs(G)^\wedge$, we first show that
for each $u \in \go$, we have $l^u$; an irreducible representation
in $\cs(G)^\wedge$. We will do this by taking the trivial
representation of the isotropy group $\cs$-algebra $\cs(G_u^u)$ and
induce that to a representation of $\cs(G)$.   To deal with
potentially discontinuous isotropy, we actually need a bit more. For
each $u \in \go$, given a subgroup $A$ of $G^u_u$ and a
representation $\pi$ of $A$, we define a representation of $\cs(G)$
and denote it $\Ind(u,A,\pi)$.  We will then let $l^u = \Ind(u,
G^u_u, 1)$ where $1$ is the trivial representation of $G^u_u$. Once
we have defined $l^u$, we will show that the map from $u \mapsto
l^u$  is indeed a well-defined, continuous, and injective map from
$\go/G$ to $\cs(G)^\wedge$ .

\section{Preliminaries}

A groupoid $G$ is a small category in which every morphism is
invertible.   Assume that $G$ is second countable, locally compact,
Hausdorff and has a fixed left Haar system,
$\{\lambda^u\}_{u\in\go}$. We define maps $r$ and $s$ from $G$ to
$G$ by $r(x) = xx^{-1}$ and $s(x) = x^{-1}x$.  These are the maps
Renault calls $r$ and $d$ in \cite{renault}. The common image of $r$
and $s$ is called the unit space which we denote $\go$.

Now consider the vector space $\cc(G)$, the space of continuous
functions with compact support from $G$ to the complex numbers,
$\C.$  We can view this space as a $*$-algebra by defining
convolution and involution with the formulae:
\begin{align}
f*g(x) &= \int f(y)g(y^{-1}x) \ d \lambda^{r(x)}(y)\notag \\
&=\int f(xy)g(y^{-1}) \ d \lambda^{s(x)}(y)\notag
\end{align}
and
\begin{equation}
f^*(x) = \overline{f(x^{-1})}.\notag
\end{equation}

We define the groupoid $\cs$-algebra with the following theorem.

\begin{theorem}
For $f \in \cc(G)$, the quantity
\begin{equation}
\|f\| := \sup \{ \|\pi(f)\| \mid \pi \text{ is a
representation of } \cc(G) \}
\label{normdef}
\end{equation}
is finite and defines a $\cs$-norm on $\cc(G)$.  The completion of
$\cc(G)$ with respect
to this norm is a $\cs$-algebra, denoted $\cs(G)$.
\end{theorem}

 This theorem can be proven as a result of
Renault's Disintegration Theorem  \cite[Theorem~4.2]{renault2},
\cite[Theorem~3.23]{muhlysnotes}.  The motivating example of
a groupoid $\cs$-algebra is a transformation group $\cs$-algebra, $\cs(H,X)$,
defined in \cite{danasthesis}.

We define the map $\pi:G \rightarrow \go \times \go$ by
$\pi(x) = (r(x),s(x))$.  Using  $\pi$, we define an equivalence
relation on $\go$ and endow the set
of equivalence classes
with the quotient topology.  We call this topological space
the orbit space of $G$, denoted
$\go/G$.

The isotropy group (also called the stability group) of $G$ at a unit $u \in \go$ is defined to be
$G^u_u := \pi^{-1}(u,u)$.  We also define $G^u := r^{-1}(u)$ and $G_u := s^{-1}(u)$.

\section{Creating a Hilbert $\cs(A)$-module}
The first step in defining a continuous injection from $\go/G$ to $\cs(G)^\wedge$ is to
define an irreducible representation $l^u$ for
each unit $u \in \go$.

Let $A$ be a subgroup of $G^u_u$ for some $u \in \go$ and let $\beta$ be a
Haar measure on $A$.
When necessary to emphasize the unit, we will sometimes write $A_u$
for $A$.
Our goal is to get an induced representation of $\cs(G)$ from a
representation
of $\cs(A)$. We will do this by showing that $\cc(G_u)$ is a Hilbert
$\cs(A)$-module on
which $\cc(G)$ acts as adjointable operators. Then we will use
\cite[Proposition~2.66]{rw-book} to achieve the desired result.

Before we get started, we need the following two technical lemmas.
\begin{remark}
We need \lemref{technicallemma} to define $l^u$, and we
need \lemref{technicallemma2}
to show that $l^u$ is irreducible.
\end{remark}

\begin{lemma}
\label{technicallemma}
Suppose that $u \in \go$ and that $A \subset G_u^u$ the isotropy
group at
$u \in \go$, and that
$\beta$ is a Haar measure on $A$.

\begin{enumerate}

\item  There is a non-negative, continuous function $b$ on
$G_u$, such that for any compact subset $K \subset G_u$,  the
support of
$b$ and $KA$ have compact intersection, and

\begin{equation}
\label{b=1}
\int_{A}b(\gamma a) d\beta(a) =1
\end{equation}

for all $\gamma \in G_u$.

\item  The formula
\begin{equation}
 Q (f)(\dot{\gamma}) = \int_{A} f(\gamma a) d\beta(a)
\end{equation}
defines a surjection from $\cc(G)$ onto $\cc(G_u/A)$.
\end{enumerate}
\end{lemma}

\begin{proof}
1.  The existence of a function $b'$ satisfying all of the
conditions of 1 except for \eqref{b=1} comes from Lemma~1 on page
96 of \cite{dixmier}.  Thus if we define
\begin{equation}
    b(\gamma ) = \dfrac{b'(\gamma )}{\tilde{b'}(\dot{\gamma})} \notag
\end{equation}
where $\tilde {b'}(\dot{\gamma}) = \int_{A} b'(\gamma a) \
d\beta(a)$.
Notice that $\tilde{b'}$ is everywhere strictly greater
than zero because \cite{dixmier} tells us that
$b'\geq 0$, $b'$ is continuous,
and $b'$ is not equal to zero on any entire equivalence class.

2.  To see that $Q(f) \in \cc(G_u/A)$ for $f \in \cc(G)$, we must
show that $Q(f)$ is continuous and has compact support.

First, we will show $Q(f)$ is continuous.  It suffices to show that
if
$\gamma_n$ converges to $\gamma$ in $G_u$ then
$Q(f)(\dot{\gamma_n})$ converges to $Q(f)(\dot{\gamma})$.
Consider the function $F_\gamma \in \cc(A)$ defined by
\begin{equation}
F_\gamma (a) := f(\gamma a).
\notag
\end{equation}
We claim that $F_{\gamma_n}$ converges to $F_{\gamma}$ in the
inductive limit topology on $\cc(A)$. It suffices to show that there
is a compact set $K \subset A$ such that $\supp F_{\gamma_n}$ is
eventually contained in $K$ and that $F_{\gamma_n}$ converges
uniformly to $F_\gamma$.

Since $\{\gamma_n\}$ converges to $\gamma$, we can pass to
a subsequence, relabel, and assume that $\{\gamma_n\} \in
\bar{K}$ for some compact set $\bar{K} \subset G_u$.
Notice that for every $n$,
\begin{equation}
\supp F_{\gamma_n}= \{a \in A \mid \gamma_n a \in \supp f\}
\subset \bar{K}^{-1}(\supp f) \cap A=K.
\notag
\end{equation}
Also note that $K$ is compact.

Now we must show that $\{F_{\gamma_n}\}$ converges to
$F_\gamma$ uniformly.  Suppose the contrary.  Let $\epsilon$ be given.
There exists a sequence $\{a_n\} \subset A$ so that
for every $n$,
\begin{equation}
|F_{\gamma_n}(a_n) - F_\gamma(a_n)| \geq \epsilon.
\label{bob}
\end{equation}
In order for (\ref{bob}) to be
greater than zero,  $\{a_n\} \subset K$.  Since $K$
is compact, we can
assume that $\{a_n\}$ converges to some $a \in K$.
Notice that $\{\gamma_n a_n\}$ converges to $\gamma a$ and
$\{\gamma a_n\}$ also converges to $\gamma a$.  Thus we can
find $N$ such that if $n \geq N$ the following holds:
\begin{align}
|F_{\gamma_n}(a_n) - F_\gamma(a_n)|&=
|f(\gamma_n a_n) - f(\gamma a_n)|\notag\\
&\leq|f(\gamma_n a_n) -  f(\gamma a )|+
|f(\gamma a)- f(\gamma a_n)|\notag\\
&\leq \epsilon\notag
\end{align}
which is a contradiction.
Therefore, $\{F_{\gamma_n}\}$ converges to
$F_\gamma$ uniformly which means
they converge in the inductive limit topology as claimed.
It follows from the continuity of
integration that $\{Q(f)(\dot{\gamma_n})\}$ converges to
$Q(f)(\dot{\gamma})$.

To see that $\supp Q(f)$ is compact, let
$\dot{\gamma} \in \supp Q(f)$.  Then $f(\gamma a) \neq 0$ for some $a$.
Thus $\gamma a \in \supp f$ which means
$\gamma a \in J= \supp f \cap G_u$, a compact subset of $G_u$.
This means that $\dot{\gamma} \in JA$, a compact subset of $G_u/A$.

We now must show that $Q$ is surjective.  Let $h \in \cc(G_u/A)$.  We can
extend $h$ and view it as a function on $G$ that is constant on orbits and
zero elsewhere.  Now consider the function $hb$ where $b$ is the
function from part 1.  Notice that
\begin{align}
\supp {(hb)} &= \supp h \cap \supp b\notag\\
&=KA \cap \supp b\notag
\end{align}
which is compact by part 1.  Now we have
\begin{align}
    Q(hb)(\dot{\gamma}) &= \int_{A} hb(\gamma a) \ d \beta(a)\notag\\
    &= \int_{A} h(\gamma a)\ d\beta(a) \int_{A} b(\gamma a)\ d \beta
    (a)\notag
\end{align}
Since $\int_{A} b(\gamma a)\ d\beta(a) = 1$ for $\dot{\gamma}$, we
have
\begin{equation}
Q(hb)(\dot{\gamma})= \int_{A} h(\gamma a)\ d\beta(a).\notag \notag
\end{equation}
But $h$ is constant on orbits, so $\int_{A}h(\gamma a)\
d\beta(a)=h(\dot{\gamma})$.
\end{proof}

Let $b'$ be the function defined in the proof of
\lemref{technicallemma}.  Now define
\begin{equation}
\rho(\gamma) = \int_{A} b'(\gamma a) \Delta(a)^{-1}
\ d\beta(a)
\label{rhodef}
\end{equation}
where $\gamma \in G_u$ and $\Delta$ is the modular function on $A$.
Notice that the  image of
$\rho(\gamma) >0$ for all $\gamma$ because
the modular function is greater than zero, and $b'$ is non-negaitve
and
non-zero on each equivalence class.  Also
\begin{align}
    \rho(\gamma c) &= \int_{A} b'(\gamma c a) \Delta(a)^{-1} \
    d\beta(a)
    \notag\\
    &=\int_{A}b'(\gamma a)\Delta(c) \Delta(a)^{-1} \ d\beta(a)
    \notag\\
    &=\Delta(c)\rho(\gamma).\label{eq:10}
\end{align}

\begin{lemma}
\label{technicallemma2}
There is a radon measure $\sigma$ on $G_u/A$ such that
\begin{equation}
\int_{G_u} f(\gamma)\rho(\gamma)\ d\lambda_u(\gamma)
=\int_{G_u/A} \int_{G^u_u} f(\gamma a) \
d\beta(a)d\sigma(\dot{\gamma})
\end{equation}
for all $f \in \cc(G_u)$.
\end{lemma}

\begin{proof}
The Riesz Representation Theorem from \cite[Theorem~7.2]{folland}
tells us that it suffices to show that the equation
\begin{equation}
\sigma(Q(f)) = \int_{G_u} f(\gamma)\rho(\gamma)\
d\lambda_u(\gamma)
\notag
\end{equation}
is a well-defined, positive, linear functional on $\cc(G_u/A)$.

It is clear that $\sigma$ is positive because $\rho$ is positive.
It is also
clear that $\sigma$ is linear.  We must show
that $\sigma$  is well-defined.  For this, it suffices to show that if
\begin{equation}
Q(f) = \int_{A} f(\gamma a)\ d\beta(a) = 0\notag
\end{equation}
for all $\gamma \in G_u$, then
\begin{equation}
    \int_{G_u} f(\gamma)\rho(\gamma) \ d\lambda_u(\gamma)=0.
\notag
\end{equation}
Let $h \in \cc(G)$.
Notice that
\begin{align}
\int_{A} h*f\ d\beta(a) &= \int_{A}\int_G h(\gamma)
f(\gamma^{-1}a) \ d\lambda^u(\gamma)d\beta(a)\label{kipper}\\
&=\int_G h(\gamma^{-1}) \int_{A}f(\gamma a)\
d\beta(a)d\lambda_u(\gamma)=0.\notag
\end{align}

Now let $\tau \in \cc(G_u/A)$ be such that $\tau(\dot{\gamma})=1$
for $\gamma \in \supp(f)$.  We can define such a function and still
require it have compact support because the image of the compact
set, $\supp f$, under the quotient map is compact. Thus the Tietze
Extension Theorem tells us we can extend this to a function in
$\cc(G_u/A)$.

Now consider
\begin{equation}
\tilde{h}(\gamma) = \tau(\dot{\gamma})b(\gamma)\rho(\gamma).
\notag
\end{equation}
We claim that this function is an element of $\cc(G_u)$.
To see this, notice that if $\gamma \in \supp \tilde{h}$ then
\begin{equation}
\gamma \in (\supp \tau \cap \supp\rho) \cap \supp b.
\notag
\end{equation}
We know that $\supp\tau$ is compact.  Thus $\supp \tilde{h}$ is also compact
because $\supp b$ intersected with a compact set is compact by
Lemma~\ref{technicallemma}.

Also, if $\gamma \in \supp f$ then
\begin{align}
\tilde{h}(\gamma a) &= \tau(\dot{\gamma})b(\gamma a) \rho(\gamma a)
\notag\\
&= b(\gamma a)\rho(\gamma a).
\notag
\end{align}

Now define
\begin{equation}
h(\gamma) := \tilde{h} (\gamma^{-1})
\notag
\end{equation}
and notice that $h \in \cc(G)$.
Thus (\ref{kipper}) implies
\begin{align}
0&= \int_{A}h*f(a)\ d\beta(a)\notag\\
&=\int_{A}\int_G h(\gamma)f(\gamma^{-1} a)\
d\lambda^u(\gamma)d\beta(a)
\notag\\
&=\int_{A}\int_G h(a\gamma)f(\gamma^{-1})\
d\lambda^u(\gamma)d\beta(a)
\notag\\
&=\int_G\int_{A} h(a\gamma^{-1})f(\gamma)\ d\beta(a)
d\lambda_u(\gamma)
\notag\\
&=\int_G\int_{A} \tilde{h}(\gamma a^{-1})
f(\gamma)\ d\beta(a)d\lambda_u(\gamma)
\notag\\
&=\int_G\int_{A} b(\gamma a^{-1}) \rho(\gamma a^{-1})\ d\beta(a)
f(\gamma)\ d\lambda_u(\gamma),
\label{lala}
\end{align}
where (\ref{lala}) follows from the fact that $\gamma \in \supp f$.
Now we get (\ref{lala}) equal to
\begin{align}
&=\int_G\int_{A} b(\gamma a^{-1})\Delta(a^{-1}) \rho(\gamma)\
d\beta(a)
f(\gamma)\ d\lambda_u(\gamma)\notag\\
&=\int_G \int_{A} b(\gamma a) \ d\beta(a) \rho(\gamma) f(\gamma)\
d\lambda_u
\notag\\
&=\int_{G_u}  f(\gamma) \rho(\gamma)\ d\lambda_u
\notag
\end{align}
as needed.
\end{proof}

Now we can begin showing that $\cc(G_u)$ is a Hilbert $\cs(A)$
module, where
$A$ is a subgroup of $G^u_u$.
Define a map
$p:\cc(G_u) \rightarrow \cc(A)$ by
\begin{equation}
p(f)(a) = \Delta(a)^{-1/2}f(a)
\notag
\end{equation}
where $\Delta$ is the modular function of $A$.  Recall that for
$\phi, \xi \in \cc(A)$,
\begin{align}
    \phi^*(a) &= \Delta(a^{-1})\overline{\phi(a^{-1})}\text{
and}\notag\\
    \xi * \phi (a) &= \int_A \xi(r) \phi(r^{-1}a) \ d
\beta(r).\notag
\end{align}
Note that $p(f^*) = p(f)^*$ where $f^*(\gamma) =
\overline{f(\gamma^{-1})}$ for $f \in \cc(G_u)$.

We know that there is a one-to-one correspondence between
unitary representations of $A$ and non-degenerate representations
of $\cs(A)$.  As is typical, we will not distinguish between these
representations with our notation.  Whether we mean the
representation
of $\cs(A)$ or of $A$ should be clear from context.

Notice that $\cc(G_u)$ is a $\cc(G)-\cc(A)$ bimodule.
The space
$\cc(A)$ acts on the right of $\cc(G_u)$ by the map
$(g,\phi)
\mapsto g
\cdot \phi:\cc(G_u) \times
\cc(A) \rightarrow \cc(G_u)$ defined by
\begin{equation}
g \cdot \phi(\gamma) =
\int_A \Delta(a)^{-1/2} f(\gamma a) \phi(a^{-1}) \ d\beta(a)\notag
\end{equation}
and $\cc(G)$ acts on the left of $\cc(G_u)$ simply by convolution
(in
$G$):
\begin{equation}
f \cdot g(\gamma) = g*f(\gamma)
\notag
\end{equation}
for $g \in \cc(G_u)$ and  $f \in \cc(G)$.
It is routine to verify that
\begin{align}
    f\cdot (\psi * \phi) &= (f \cdot \psi) \cdot \phi \text{
and}\notag\\
    (f*g) \cdot g &= f \cdot (h \cdot g).\notag
\end{align}

\begin{lemma}
\label{hibertmodule}
The formula
\begin{equation}
\langle f,g \rangle_A = p(f^**g) \notag
\end{equation}
defines a $\cs(A)$-valued sesquilinear form on $\cc(G_u)$ such
that for all $f,g \in \cc(G_u)$ and $\phi \in \cc(A)$ the following
conditions hold:
\begin{enumerate}
\item $\langle f,g \rangle_A^* = \langle g,f \rangle_A$
\item $\langle f, g \cdot \phi \rangle_A = \langle f,g \rangle_A
* \phi$
\item $\langle f,f \rangle_A $ is a positive element of $\cs(A)$.

\end{enumerate}
\end{lemma}

\begin{proof}
It is straightforward to check that this map is sesquilinear.

1. Notice that
\begin{align}
\langle f,g \rangle_A^*(a) &= (p(f^* * g))^*(a)\notag\\
&=\Delta(a^{-1})\overline{p(f^* * g)(a^{-1})}\notag\\
&=\Delta(a^{-1})\Delta(a^{-1})^{-1/2}\overline{(f^**g)(a^{-1})}\notag\\
&=\Delta(a)^{-1/2}(f^**g)^*(a)\notag\\
&= \Delta(a)^{-1/2}(g^**f)(a)\notag\\
&=\langle g,f \rangle _A(a).\notag
\end{align}

2.  We see that:
\begin{align}
\langle f,g \rangle_A * \phi  (a) &= \int_A \langle f,g \rangle
_A(r)
\phi(r^{-1}a)\ d\beta(r)\notag\\
&= \int_A \Delta(r)^{-1/2} f^* * g(r) \phi (r^{-1}a) \
d\beta(r)\notag\\
&=  \int_A \Delta(r)^{-1/2} \int_G f^*(x)  g(x^{-1}r)  \phi
(r^{-1}a)\
d\lambda^u
(x)d\beta(r).\notag
\end{align}
Working on the right hand side of the equation we want to prove, we
get:
\begin{align}
\langle f,g \cdot \phi \rangle _A (a) &= \Delta(a)^{-1/2} (f^* *
(g \cdot
\phi)(a)\notag\\
&= \Delta(a)^{-1/2} \int_G f^*(x) (g\cdot \phi)(x^{-1}a)\
d\lambda^u(x)\notag\\
&= \Delta(a)^{-1/2} \int_G f^*(x) \int_A \Delta(r)^{-1/2}
g(x^{-1}ar)
\phi(r^{-1})\ d\lambda^u(x)\label{dipsy}
\end{align}
multiplying $r$ on the left by $a^{-1}$ we see that
equation~(\ref{dipsy}) is equal to:
\begin{align}
&= \Delta(a)^{-1/2} \int_G f^*(x) \int_A \Delta(a^{-1}r)^{-1/2}
g(x^{-1}r)
\phi((ar)^{-1})\ d\lambda^u(x)\notag\\
&= \Delta(a)^{-1/2} \int_G f^*(x) \int_A
\Delta(a)^{1/2}\Delta(r)^{-1/2} g(x^{-1}r)
\phi((r^{-1}a)\ d\lambda^u(x)\notag\\
&=  \int_A \Delta(r)^{-1/2} \int_G f^*(x)  g(x^{-1}r)  \phi
(r^{-1}a)\
d\lambda^u
(x)d\beta(r).\notag
\end{align}

Thus,
\begin{equation}
\langle f,g \cdot \phi \rangle_A = \langle f,g\rangle_A * \phi. \notag
\end{equation}

3.  To show positivity,
it suffices to show that
\begin{equation}
(\pi(\langle f ,f \rangle_A)h \mid h ) \geq 0
\notag
\end{equation}
for any irreducible representation
$\pi \in \cs(A)^\wedge$.

Let $\pi \in \cs(A)^\wedge$.  Define a
Radon measure $\mu$ on $A$  by
\begin{align}
    \int_A \phi(a) \ d\mu(a)  &= \int_A \phi(a) (\pi(a) \xi \mid \xi)\
    d\beta(a)\notag\\
    &= (\pi(\phi)\xi \mid \xi)\notag
\end{align}
for $\phi \in \cc(A)$.
Notice that $\phi \geq 0$ if and only if
$\int_A\phi(a)\ d\mu(a) \geq 0$ for every $\pi \in \hat{A}$
and $\xi \in \H_\pi$.

Now consider
\begin{align}
(\pi(\langle g,f \rangle_A)h \mid h ) &= \int_A \langle g,f
\rangle_A(a)((\pi(a)h \mid h ) \ da \notag\\
&= \int_A  \langle g,f
\rangle_A (a) \ d\mu(a) \notag \\
&= \int_A \Delta(a)^{-1/2} g^* *f(a)\ d\mu(a)\notag\\
&=\int_A \Delta(a)^{-1/2} \int_G g^*(\gamma)f(\gamma^{-1}a) \
d\lambda^u(\gamma)d\mu(a)\notag\\
&= \int_A \Delta(a)^{-1/2}
\int_G \overline{g(\gamma^{-1})}f(\gamma^{-1}a)
\ d\lambda^u(\gamma)d\mu(a).\notag
\end{align}
We can change $\gamma^{-1}$ to $\gamma$ by switching the
measure from
$\lambda^u$ to $\lambda_u$ and get:
\begin{align}
&= \int_A \Delta(a)^{-1/2}
\int_G \overline{g(\gamma)}f(\gamma a)
\ d\lambda_u(\gamma)d\mu(a)\notag\\
&= \int_A \Delta(a)^{-1/2}
\int_G \overline{g(\gamma)}f(\gamma a) \int_A b(\gamma t ) \
d\beta(t) \
\ d\lambda_u(\gamma)d\mu(a)\notag
\end{align}
where $b$ is from \lemref{technicallemma} so that
$\int_A b(\gamma t ) \ d\beta(t) =1$ for all $\gamma \in G_u$.
Using Fubini's theorem, we get this equal to
\begin{align}
&= \int_A \Delta(a)^{-1/2} \int_A \int_G
\overline{g(\gamma)}f(\gamma a) b(\gamma t ) \
d\lambda_u(\gamma)d\beta(t)d\mu(a)\notag\\
&= \int_A \Delta(a)^{-1/2} \int_A \int_G
\overline{g(\gamma^{-1})}f(\gamma^{-1} a) b(\gamma^{-1} t ) \
d\lambda^u(\gamma)d\beta(t)d\mu(a).\notag
\end{align}
Since $t \in A \subset G^u_u$, we can multiply $\gamma$ on the left
by
$t^{-1}$ to see that this is equal to
\begin{align}
&= \int_A \Delta(a)^{-1/2} \int_A \int_G
\overline{g(\gamma^{-1}t^{-1})}f(\gamma^{-1}t^{-1} a)
b(\gamma^{-1}) \
d\lambda^u(\gamma)d\beta(t)d\mu(a)\notag\\
&=\int_A \Delta(a)^{-1/2} \int_A \int_G
\overline{g(\gamma t^{-1})}f(\gamma t^{-1} a) b(\gamma) \
d\lambda_u(\gamma)d\beta(t)d\mu(a)\notag\\
&= \int_G b(\gamma) \int_A \int_A \overline{g(\gamma
t^{-1})}\Delta(a)^{-1/2}f(\gamma t^{-1} a) \
d\beta(t)d\mu(a)d\lambda_u(\gamma)\notag\\
&= \int_G b(\gamma) \int_A
\int_A\Delta(a)^{-1/2} \overline{g(\gamma t^{-1})}f(\gamma t^{-1}
a)
\ d\beta(t)d\mu(a)d\lambda_u(\gamma).\label{eq10}
\end{align}
Recall that $\Delta(t^{-1}t)=\Delta(e) = 1$. So we have
equation~\ref{eq10}
\begin{align}
&= \int_G b(\gamma) \int_A
\int_A\Delta(t^{-1}t)\Delta(a)^{-1/2}
\overline{g(\gamma t^{-1})}f(\gamma t^{-1} a)
\ d\beta(t)d\mu(a)d\lambda_u(\gamma)\notag\\
&= \int_G b(\gamma) \int_A
\int_A\Delta(t^{-1})\Delta(t)^{1/2}
\overline{g(\gamma t^{-1})}\Delta(t^{-1})^{-1/2}\Delta(a)^{-1/2}
f(\gamma t^{-1} a)
\ d\beta(t)d\mu(a)d\lambda_u(\gamma)\notag\\
&= \int_G b(\gamma) \int_A
\int_A\Delta(t^{-1})\Delta(t)^{1/2}
\overline{g(\gamma t^{-1})}\Delta(t^{-1}a)^{-1/2}f(\gamma t^{-1} a)
\ d\beta(t)d\mu(a)d\lambda_u(\gamma).\label{eq11}
\end{align}
Now define $g_\gamma, f_\gamma \in \cc(A)$ by
\begin{align}
g_\gamma(a) &= \Delta(a)^{-1/2}g(\gamma a)\text{ and}\notag\\
f_\gamma(a) &= \Delta(a)^{-1/2}f(\gamma a).\notag
\end{align}
Also notice that
\begin{align}
g_\gamma^*(t)
&= \Delta(t^{-1}) \overline {g_\gamma(t^{-1})}
\notag\\
&= \Delta(t^{-1})
\Delta(t)^{1/2} \overline{g(\gamma t^{-1})}.\notag
\end{align}
where we are viewing $p(g_\gamma)$ as an element of $\cs(A)$.

Thus, equation~(\ref{eq11}) is equal to
\begin{align}
&= \int_G b(\gamma) \int_A \int_A g_\gamma^*(t)
f_\gamma(t^{-1} a) \ d\beta(t)
d\mu(a) d\lambda_u(\gamma)\notag\\
&= \int_G b(\gamma) \int_A g_\gamma^* * f_\gamma (a) \ d\mu(a)
d\lambda_u(\gamma).\label{eq12}
\end{align}
Using the definition of $\mu$, this equals
\begin{align}
&= \int_G b(\gamma)(\pi(g_\gamma^* * f_\gamma )h \mid h) \
d\lambda_u(\gamma)
\notag\\
&= \int_G b(\gamma)(\pi(f_\gamma)h \mid \pi(g_\gamma) h) \
d\lambda_u(\gamma).
\notag\
\end{align}
Replacing $g$ by $f$  we get
\begin{align}
(\pi(\langle f, f \rangle_A)h\mid h)
&= \int_G b(\gamma)(\pi(f_\gamma)h \mid
\pi(f_\gamma) h) \ d\lambda_u(\gamma)\notag\\
&= \int_G b(\gamma)\Delta(h)^{-1/2} (\pi(f_\gamma)h \mid
\pi(f_\gamma) h) \ d\lambda_u(\gamma) \geq 0\notag
\end{align}
for all $\pi \in \hat{A}$ and $\xi \in \H_\pi$.
This completes the proof.
\end{proof}

With this lemma, we see from \cite[Lemma~2.16]{rw-book} that
$\cc(G_u)$ is
a pre-Hilbert $\cs(A)$-module.  We now want to show that $\cc(G)$ acts
on the right
of this module as adjointable operators.  In order to do that, we
need \lemref{scooby} below.  The proof to this lemma requires we use Renault's
Disintegration theorem.  The corresponding lemma in \cite[Proposition~2.2]{ctstraceIII} does not
require this difficult theorem.  However, in removing the assumption the
isotropy groups be abelian, the modular function that is introduced causes
the argument in \cite[Propsotion~2.2]{ctstraceIII} to be invalid.

The specific portion of Renault's Disintegration Theorem, which is found in
\cite[Thoerem~3.23]{muhlysnotes}, is that a representation of $\cc(G)$ on a
pre-Hilbert space $\H_0$ extends uniquely to a bounded operator of $\cc(G)$ on
the Hilbert space completion of $\H_0$, $\H$.

\begin{lemma}
\label{scooby}
The inequality
\begin{equation}
\langle f*g , f*g \rangle_A \leq \|f\|_{\cs(G)} \langle g,g \rangle
_A,
\label{part2}
\end{equation}
holds for all $f \in \cc(G_u)$ and $g \in \cc(G)$.
\end{lemma}

\begin{proof}

We know from the previous lemma that $\langle f,f,\rangle_A \geq 0$
in $\cs(A)$.  Thus, given any state
$\rho$ in $\cs(A)$,
\begin{equation}
(f | g)_\rho := \rho(\langle g , f \rangle_A)\notag
\end{equation}
defines a pre-inner product on $\cc(G)$.  The quotient $\H_0 = \cc(G)/N$ where
\begin{equation}
N:= \{ f \in \cc(G_u) \mid \rho(\langle f,f \rangle_A ) =0\}\notag
\end{equation}
is a dense subspace of the Hilbert space completion of $\H_0$ with respect to
$( \cdot|\cdot)_\rho$.  Now define a a homomorphism $L$ from $\cc(G)$ into
the algebra of linear transformations on $\H_0$ by  $L(f)g = f*g$.
We will show that $L$  satisfies conditions (a)--(c) of
\cite[Theorem~3.32]{muhlysnotes}.
\begin{enumerate}
\item[(a)]  The span of
\begin{equation}
\{L(f)g= f*g \mid f \in \cc(G), g\in \cc(G_u)\}\notag
\end{equation}
is dense in $\cc(G_u)$.  This follows from the existence of an
approximate identity \cite[Corollary~2.11]{equivalence}.
\item[(b)]  For each $g,h \in \cc(G_u)$, the functional
$L_{g,h}:\cc(G) \rightarrow \C$ defined by the equation
\begin{equation}
L_{g,h}(f) = (g|L(f)h)_\rho=(g|f*h)_\rho
\notag
\end{equation}
is continuous with respect to the inductive limit topology on $\cc(G)$
because the inner product is just a compositions of continuous functions.
Specifically, the function $p$, convolution, the
modular function, and $\rho$ are all continuous.

\item[(c)]  For all $f \in \cc(G)$, and for all $g,h \in \cc(G_u)$,
\begin{align}
(g|L(f^*)h) &= (g|f^**h)\notag\\
&=(f*g|h)\notag\\
&=(L(f)g|h)
\notag
\end{align}
as needed.
\end{enumerate}

Invoking \cite[Theroem~3.23]{muhlysnotes}, we know that
$L$ is a representation of $\cc(G)$ on $\H$.  Now we can use
\cite[3.3.7]{murphy} to get
\begin{align}
(f*g | f*g)_\rho &\leq \|f\|_{\cs(G)} (g|g)_\rho\text{ or}\notag\\
\rho(\langle f*g,f*g \rangle_A) &\leq \|f\|_{\cs(G)} \rho(\langle g,g \rangle_A)
\notag
\end{align}
for all states $\rho$.  Thus
\begin{equation}
\langle f*g,f*g \rangle_A \leq \|f\|_{\cs(G)} \langle g,g \rangle_A
\notag
\end{equation}
in $\cs(A)$.
\end{proof}

\begin{lemma}
The space
$\cc(G)$ acts on the right
of the pre-Hilbert $\cs(A)$-module $\cc(G_u)$ as adjointable
operators.
\end{lemma}

\begin{proof}
It suffices to show that
\begin{equation}
\langle f * h, g \rangle_A = \langle f, h^* * g \rangle_A
\notag
\end{equation}
for $f \in \cc(G)$ and $g,h \in \cc(G_u)$.
This equation is fairly obvious.  Notice that
\begin{align}
\langle h*f,g \rangle_A &= p((f*h)^* * g)\notag\\
&= p(h^* * f^* * g)\notag\\
&= p(h^* * (f^* * g)) \notag\\
&= \langle h, f^* * g \rangle_A \notag
\end{align}
\end{proof}

Now given a representation,
$\pi$ of $\cs(A)$ on $\H_\pi$
\begin{equation}
\langle f \otimes \xi, g \otimes \eta \rangle =
(\pi(\langle g,f \rangle_A)\xi\mid
\eta)_{\H_\pi}
\label{innerproduct}
\end{equation}
defines a pre-inner product on $\cc(G_u) \otimes \H_\pi$ on which
$\cc(G)$
acts by bounded operators via convolution by
\cite[Proposition~2.66]{rw-book}.


\section{Defining $l^u$}

\begin{definition}
Suppose that for $u \in \go$, $A$ is a subgroup of $G^u_u$, and  $\pi$
is a
representation of $A$ on $\H_\pi$.  Then we define the
representation
$\Ind(u,A,\pi)$ on the Hilbert space completion of $\cc(G_u) \otimes
\H_\pi$
with respect to the inner product determined in \ref{innerproduct}
according to \cite[Proposition~2.66]{rw-book} by:
\begin{equation}
\Ind(u,A,\pi)(f)[g \otimes \xi]  := [f*g\otimes\xi].
\end{equation}
\end{definition}

Now we will define $l^u$.

\begin{equation}
l^u := \Ind(u, G^u_u, 1)
\notag
\end{equation}
where $1$ is the trivial representation of $G^u_u$.
We will identify the Hilbert space that $l^u$ acts upon as the
Hilbert
space completion of $\cc(G_u)$ with respect to the inner product
\begin{equation}
(f|g )_u = \int _{G^u_u} g^* * f(a)\Delta(a)^{-1/2} \ d \beta(a).
\notag
\end{equation}
Thus
\begin{equation}
l^u(f)g = f*g.
\notag
\end{equation}

\begin{lemma}
\label{irreducible}
Each $l^u$ is an irreducible representation.
\end{lemma}

\begin{proof}
Following the proof of \cite[Lemma~2.5]{ctstraceIII}, we will show
that $l^u$
is irreducible, by showing that $l^u$ is equivalent  to a
representation we will call $R^u$ and that $R^u$ is irreducible.

Recall that $l^u$ is a representation on $\H_u$ which is the
Hilbert
space
completion of $\cc(G_u)$ with respect to the inner product
\begin{equation}
( f|g)_u=\int_{A_u} g^* * f(a)\Delta(a)^{-1/2} \ d \beta ^u (a).
\notag
\end{equation}
Let $\V_u = L^2(G_u/A, \sigma^u)$ where $\sigma^u$ is the the measure found in
\lemref{technicallemma2}.  Recall that $\V_u$ is a Hilbert space
with respect to the inner product
\begin{equation}
\langle \xi, \zeta \rangle_u = \int_{G_u/A_u}
\xi(\gamma)\overline{\zeta(\gamma)} \ d \sigma^u
(\dot{\gamma}).\notag
\end{equation}
Consider the map $U^u:\H_u \rightarrow \V_u$ defined by
\begin{equation}
U^u(g)(\gamma) = \int_{A_u} g(\gamma a )\rho(\gamma a)^{-1/2}
 \ d\beta^u (a)\notag
\end{equation}
where $\rho$ is defined in Equation (\ref{rhodef}).
Notice that
\begin{align}
(f|g)_{u} &= \int_{A_u} g^* *f(a) \Delta(a)^{-1/2} \ d\beta^u(a)
\notag \\
&= \int_{A_u}\int_G g^*(a\gamma)f(\gamma^{-1})\
d\lambda^u(\gamma)\Delta(a)^{-1/2} d
\beta^u(a) \notag\\
&= \int_{A_u}\int_G
\overline{g(\gamma^{-1}a^{-1})}f(\gamma^{-1})\
d\lambda^u(\gamma)\Delta(a)^{-1/2} d
\beta^u(a) \notag\\
&= \int_{A_u}\int_G \overline{g(\gamma a^{-1})}f(\gamma)\
d\lambda_u(\gamma)\Delta(a)^{-1/2} d
\beta^u(a). \label{eq:101}
\end{align}
Now we can map $a \mapsto a^{-1}$ and introduce
$\rho(\gamma)^{-1}\rho(\gamma)=1$ to
see that (\ref{eq:101}) is equal to
\begin{align}
&= \int_{A_u}\int_G \overline{g(\gamma a)}f(\gamma)
\rho(\gamma)^{-1}\rho(\gamma) \
d\lambda_u(\gamma)\Delta(a)^{-1/2} d
\beta^u(a) \notag\\
&=\int_{A_u}\int_{G/A_u} \int_{A_u} \overline{g(\gamma b
a)}f(\gamma b)
\rho(\gamma b)^{-1}  \ d\beta(b) d\sigma(\dot{\gamma})
\Delta(a)^{-1/2} d
\beta^u(a) .\label{eq:102}
\end{align}
Multiplying $a$ on the left by $b^{-1}$, using
Fubini's Theorem, and
equation (\ref{eq:10})  we
get (\ref{eq:102}) equal to
\begin{align}
&=\int_{G/A_u}\int_{A_u} \int_{A_u} \overline{g(\gamma
a)}f(\gamma b)
\rho(\gamma )^{-1}\Delta(b)^{-1}\Delta(b^{-1}a)^{-1/2}  \ d\beta(a)
d\beta(b)d\sigma(\dot{\gamma}) \notag\\
&=\int_{G/A_u}(\int_{A_u} \overline{g(\gamma  a)}\Delta(a)^{-1/2}
\rho(\gamma)^{-1/2} \ d\beta(a)) (\int_{A_u}f(\gamma b)
\rho(\gamma )^{-1/2}\Delta(b)^{-1/2}  \ d\beta(b) )\notag\\
&= \int_{G/A_u} \overline{U_u(g)(\dot{\gamma})}
U_u(f)(\dot{\gamma}) \ d\sigma(\dot{\gamma})
=\langle U_u(f),U_u(g)\rangle_u.\notag
\end{align}
This means that $U^u$ is an isometry.  To see that $U^u$ is
surjective, let $\xi \in L^2(G_u/A,\sigma^u)$.  We define a
function, $\bar{\xi} \in \cc(G_u)$ with the formula
\begin{equation}
\bar{\xi}(\gamma) = \rho(\gamma)^{1/2}b(\gamma)
\xi(\dot{\gamma})
\notag
\end{equation}
where $b$ is the function defined in \lemref{technicallemma}.
Now we see that $U^u(\bar\xi(\gamma)) \xi$ thus $U^u$ is surjective making it a unitary
operator.  We use $U^u$ to intertwine $l^u$ with the
operator $T^u$ defined below.
\begin{equation}
T^u(f)\xi(\gamma) = \int_G f(\gamma\eta)\xi(\eta^{-1}) \ d
\lambda^u(\eta)\notag
\end{equation}
where $\xi \in \V_u$.  To see this, notice
\begin{align}
T^u(f)U^ug(\gamma) &= \int_G f(\gamma\eta)\int_{A_u}g(\eta^{-1}a)
\rho(\gamma a)^{-1/2} \
d\beta^u(a)d(\lambda^u(\eta)\notag\\
&= \int_G\int_{A_u}f(\gamma\eta)g(\eta^{-1}a)\rho(\gamma
a)^{-1/2} \
d\beta^u(a)d\lambda^u(\eta)\notag
\end{align}
for $f \in \cc(G)$ and $g \in \cc(G_u)$ and $\gamma \in G_u$.  On the other hand,
\begin{align}
U^ul^u(f)g(\gamma) &= U^u(f*g)(\gamma)\notag\\
&=\int_{A_u} f*g(\gamma a) \rho(\gamma a)^{-1/2} \
d\beta^u(a)\notag\\
&=\int_{A_u}\int_G f(\gamma a \eta)g(\eta^{-1})\rho(\gamma
a)^{-1/2} \
d\lambda^u(\eta)d\beta^u(a).\notag
\end{align}
Using Fubini's theorem and multipling $\eta$ on the left by
$a^{-1}$
we get
\begin{equation}
U^ul^u(f)g(\gamma)=\int_G\int_{A_u}f(\gamma\eta)
g(\eta^{-1}a)\rho(\gamma a)^{-1/2} \
d\beta^u(a)d\lambda^u(\eta)\notag
\end{equation}
as needed.

We know  from \cite[Corollary~3.2.1]{arveson}
that $G_u$ is an analytic
Borel space.  Therefore from \cite[Theorem~3.4.1]{arveson}, there
is a
Borel cross  section $c$
for the natural map from $G_u$ to $G_u/A_u$.  Now define an
operator,
$W^u:\V_u
\rightarrow L^2(G_u/A_u,
\sigma^u)$ by
$W^u(f)(\dot\gamma)=f(c(\dot\gamma))$.  We claim that
$W^u$ is a unitary operator.  To see this, notice that
$W^*(\xi)(\gamma)=\xi(\dot\gamma)$ for
$\xi
\in L^2(G_u/A_u,\sigma^u)$.

The point here is that now we have a unitary to intertwine $T^u$
with the
representation, $M^u$ on $L^2(G_u/A_u,\sigma^u)$ given by:
\begin{equation}
M^u(f)\xi(\dot\gamma) = \int_G f(\gamma \eta^{-1})\xi(\dot\eta)
\ d
\lambda_u(\eta).\notag
\end{equation}
In showing this equivalence, we use the assumption that functions
in
$\V_u$ are constant on orbits.
Now we know $T^u$ is unitarily equivalent to $M^u$ thus  $l^u$ is unitarily
equivalent to $M^u$.

Again, since $G_u$ is
second
countable, the
 restriction of $r$ defines a Borel isomorphism of $G_u/A_u$ with
$[u]$.  We let
$\sigma_*^u = r_*(\sigma^*)$ so that
\begin{equation}
\int\phi(v) \ d\sigma_*^u(v) = \int_{G_u/A_u}\phi(r(\gamma)) \
d\sigma^u(\dot{\gamma}).\notag
\end{equation}
Then $M^u$ is equivalent to the representation $R^u$ on
$L^2([u],\sigma_*^u)$ given by
\begin{equation}
R^u(f)\xi(v)= \int_G f(\eta^{-1})\xi(r(\eta)) \
d\lambda_{v}(\eta).\notag
\end{equation}
If we let $v=r(\gamma)$ this becomes
\begin{equation}
R^u(f)\xi(r(\gamma))= \int_G f(\eta^{-1})\xi(r(\eta\gamma)) \
d\lambda_{r(\gamma)}(\eta).
\label{R^u=}
\end{equation}
Now let $N^u$ be the
representation of $\co(\go)$ on $L^2([u],\sigma_*^u)$ given by
\begin{equation}
N^u(\phi)(\xi)(v) = \phi(v)\xi(v)\notag
\end{equation}
where $v \in [u], \phi \in \co(\go)$, and $\xi \in
L^2([u],\sigma_*^u)$.

Suppose $\phi \cdot f
(\gamma) = \phi(r(\gamma))f(\gamma)$ for $ \phi \in \cc(\go) $
and $f\in
\cc(G)$. Notice that
$R^u(\phi
\cdot f) = N^u(\phi)R^u(f)$.
We will now show that $R^u$ is irreducible.
It suffices to show that
\begin{equation}
R^u(\cs(G))' = \C I.
\notag
\end{equation}

Note that $R^u(\cs(G))'$ is a von Neumann algebra.
We will consider
any projection of $\B(L^2([u],\sigma_*^u)$ that commutes with
$R^u(\cs(G))$.
If we can show that this projection is either $0$ or $1$ we are
done.
A projection is sufficient because of \cite[4.1.11]{murphy}.

Notice that any projection commuting with $R^u(\cs(G))$ must also
commute with $N^u(\co(\go))''$.
By \cite[Theorem~2.2.1]{arveson}, we know that
\begin{equation}
N^u(\co(\go))'' = L^\infty (\go, \sigma^u_*).\notag
\end{equation}
Therefore any projection commuting with $R^u(\cs(G))$ must be of
the form $N^u(\phi)$
with $\phi = 1_E$ and $E \subset [u]$ because $\supp \sigma^u_*
\in [u]$.
We will show that $E$ is either the empty set or all of $[u]$.

Since $N^u(\phi)$ commutes with $R^u(f)$ we have
\begin{align}
N^u(\phi)R^u(f)\xi(v) &= R^u(f)N^u(\phi)\xi(v)\notag\\
\phi(v) \int_G f(\eta^{-1})\xi(r(\eta))\ d\lambda_v(\eta)
&= \int_G f(\eta^{-1})\phi(r(\eta))\xi(r(\eta))\
d\lambda_v(\eta)\notag
\end{align}
for $\sigma^u_*$-almost every $v$,
$\xi \in L^2([u],\sigma^u_*$, and all $f \in
\cc(G)$.
Thus for some $v \in [u]$, $\phi(v) = \phi(r(\eta))$
$\lambda_v$-almost
every $\eta \in G_v$.
This means that $\phi$ is constant (a.e.) on $[u]$.  Thus, $E$ is
empty or all of
$[u]$ and $R^u$ is irreducible and hence $l^u$ is irreducible.
\end{proof}

We have now verified that $l^u$ is indeed an element of
$\cs(G)^\wedge$.
\section{The Continuous Injection}

Define a map $\omega: \go/G \rightarrow \cs(G)^\wedge$ so that $u
\mapsto l^u$ where $l^u(f)(g)=f*g$ as defined above.
We are being sloppy with our notation.  We actually are viewing $l^u$ as $[l^u]$, the unitary equivalence class
of $l^u$.
\begin{lemma}
The map $\omega$ defined above is well-defined.
\end{lemma}

\begin{proof}
Suppose $[u]=[v]$.   Thus there exists $\gamma \in G$ such that
$s(\gamma)=v$ and $r(\gamma)=u$.

We want to show $l^u$ is unitarily equivalent to $l^v$.
Notice that $a \mapsto \gamma a \gamma^{-1}$ is an isomorphism from
$G^v_v$
to $G^u_u$.
Also notice the map
\begin{equation}
\phi \mapsto \int_{A_v} \phi(\gamma a \gamma^{-1})\ d\beta^v(a)\notag
\end{equation}
defines a Haar measure on $G^u_u$ for all $\phi \in C_c(G^u_u)$.
This Haar measure may not equal $\beta^u$.  However, there is a
constant $z(\gamma)$ such that $z(\gamma) >0$ and
\begin{equation}\int_{A_u}\phi(a)\ d\beta^u(a) =z(\gamma) \int_{A_v}
\phi(\gamma a \gamma^{-1})\ d\beta^v(a).\notag\notag
\end{equation}

Now define $q:\H_u \rightarrow \H_v$ by $q(g)(b) = z(\gamma)^{1/2} g(b \gamma^{-1} )$
for $g \in C_c(G_u)$ and $b \in G_v$.
We claim that $q$ is a unitary.  We must be
careful in our computations in order to distinguish
the $*$ operation in $\cc(G_u)$ from that in $\cc(G^u_u)$.  For $g,h \in \cc(G_u)$, consider
\begin{align}
(q(g)|q(h))_v &= \int_{G^v_v} q(g)^**q(h)(a)(a)\ d\beta^v(a)\notag\\
&= \int_{G^v_v} \int_G q(g)^*(x)q(h)(x^{-1}a)\
d\lambda^{r(a)}(x)\Delta(a)^{-1/2}\ d\beta^v(a)\notag\\
&= \int_{G^v_v}z(\gamma) \int_G \overline{g(x^{-1}\gamma^{-1})}
h(x^{-1}a\gamma^{-1})\ d\lambda^{r(a)}(x)\Delta(a)^{-1/2}\
d\beta^v(a).\notag\\
&= \int_{G^v_v}z(\gamma) \int_G g^*(\gamma x)h(x^{-1}a
\gamma^{-1})\ d\lambda^{r(a)}(x)\Delta(a)^{-1/2}\ d\beta^v(a).\notag
\end{align}
We can multiply $x$ on the left by $\gamma^{-1}$ and get
\begin{flalign}&= \int_{G^v_v}z(\gamma)
\int_G g^*(x)h(x^{-1}\gamma a\gamma^{-1})\
d\lambda^{r(a)}(x)\Delta(a)^{-1/2}\ d\beta^v(a)\notag\\
&= \int_{G}z(\gamma) \int_{G^v_v} g^*(x)h(x^{-1}\gamma
a\gamma^{-1}) \Delta(a)^{-1/2}\ d\beta^v(a)
d\lambda^{r(a)}(x)\notag\\
&= \int_{G}z(\gamma) \int_{G^v_v} g^*(x)h(x^{-1}\gamma
a\gamma^{-1})\Delta(\gamma)^{1/2}\Delta(\gamma)^{-1/2}
\Delta(a)^{-1/2}\ d\beta^v(a)d\lambda^{r(a)}(x)\notag\\
&= \int_{G} \int_{G^u_u} g^*(x)h(x^{-1}a) \Delta(a)^{-1/2}\
d\beta^v(a)d\lambda^{r(a)}(x)\notag\\
&= \int_{G^u_u} \int_{G} g^*(x)h(x^{-1}a) \ d\lambda^{r(a)}
(x) \Delta(a)^{-1/2}\ d\beta^v(a)\notag\\
&=\int_{G^u_u}g^**h(a)\Delta(a)^{-1/2}\ d\beta^u(a)\notag\\
&=(g|h)_u.\notag
\end{flalign}This means that $q$ is isometric.
Clearly $q$ is surjective, making $q$ a unitary as claimed. Suppose
that $g \in \cc(G_u)$ and $f \in \cc(G)$, then
\begin{align}
q(l^u(f)g)b&=z(\gamma)^{1/2}l^u(f)g(b\gamma^{-1})\notag\\
&=z(\gamma)^{-1/2}f*g(b\gamma^{-1})\notag\\
&=\int_G z(\gamma)^{1/2}f(x)g(x^{-1}b\gamma^{-1})\
d\lambda^u(x)\notag\\
&=\int_G f(x)q(g)(x^{-1}b)\ d\lambda^u(x)\notag\\
&=f*q(g)(b)\notag\\&=l^v(f)(q(g))b;\notag
\end{align}
therefore, $l^u$ is unitarily equivalent to $l^v$.
\end{proof}

Before proving that our map is continuous, we need the following
lemma which is a standard result taken from
\cite[Theorem~5.9]{rieffel} and
\cite[Proposition~6.26]{rieffel}.

\begin{lemma}
\label{rief}
Suppose that $u \in \go$, that $A$ and $B$ are subgroups of $G^u_u$
with
$A \subset B$, and that $\pi$ and $\rho$ are representations of $A$.
\begin{enumerate}
\item $\Ind(u, A, \pi)$ is unitarily equivalent to
$\Ind(u,B,\Ind^B_A(\pi))$, and
\item if $\pi$ weakly contains $\rho$, the $\Ind(u,A,\pi)$ weakly
contains $\Ind(u,A,\rho)$.
\end{enumerate}
\end{lemma}

Recall that we can endow the space of subgroups of $G$ with the
compact topology described by Fell in \cite{fell}.
We will denote this space by $\Sigma$.
Thus we can view a subgroup $A$ of $G^u_u$
as an element of $\Sigma$.
Now, we can fix Haar measures $\beta^{A}$
on each subgroup $A$ in $\Sigma$ such that the map
\begin{equation}
A \mapsto \int_A f \ d \beta^{A}
\notag
\end{equation}
is continuous for each $f \in \cc(G)$.
This is a result of \cite[Lemma~1.6]{renault16}.
For the sake of simplicity, we will write $\beta^u$ for
$\beta^{G^u_u}$.

\begin{lemma}
Suppose that $G$ is a groupoid where the isotropy groups are
amenable.  Then the map  $\omega: \go/G \rightarrow \cs(G)^\wedge$
is
continuous.
\end{lemma}

\begin{proof}
It suffices to show that if $l^u(f) \neq 0$ and if $\{u_n\}$
converges
to $u$, then, eventually $l^{u_n}(f) \neq 0$.
If the assertion fails,
then it fails for some $f \in C_c(G)$.
Since $\Sigma$ is compact, we can pass to a subsequence, relabel
and assume
that $G^{u_n}_{u_n}$ converges to $C \in \Sigma$.  Note that
$C \subset G^u_u$. Passing to yet another subsequence
and relabeling, we can also assume
that $l^{u_n}(f)=0$ for all $n$.

It follows from our choice of Haar measures that, for each $g,h \in
C_c(G_u)$,
\begin{equation}
\int_{G^{u_n}_{u_n}}\Delta(a)^{-1/2} h^* * f * g(a) \ d \beta^{u_n}(a) \text{ converges to }
\int_{C} \Delta(a)^{-1/2} h^* *f*g(a) \ d \beta^C(a)
\notag
\end{equation}
which means
\begin{equation}
(l^{u_n}(f)(g)|h)_u \text{ converges to } (\Ind(u,C,1)(f)(g)|h)_u.
\notag
\end{equation}

Therefore $\Ind(u,C,1)(f) = 0$.  By \lemref{rief} part 1,
$\Ind(u,G^u_u,\Ind^{G^u_u}_C(1))(f)=0$.
 Because we assume that the isotropy groups
are amenable, from \cite[Theorem~5.1]{greenleaf}, we know that
$\Ind^{G^u_u}_C(1)$ is weakly contained in the trivial
representation of
$G^u_u$.  Thus part 2 of
\lemref{rief} implies that
$\Ind(u,G^u_u,1)(f)=l^u=0$ as desired.
\end{proof}

\begin{lemma}
The map
$\omega$ is injective.
\end{lemma}

\begin{proof}
Suppose $\omega([u]) = \omega([v])$.  That is that $[l^u] = [l^v]$
where $u,v
\in \go$.  We must show that
$[u] = [v]$.

We know there exists a unitary representation $W$ that intertwines
$l^u$ and $l^v$.  That means that $Wl^u=l^vW$.   Recall from the
proof of \lemref{irreducible} we found representations $R^u$ and $
R^v$ and representation $N^u$ and $N^v$ so that when $l^u$ is
unitarily equivalent to $l^v$ then $R^u$ is unitarily equivalent to
$R^V$ which means that $N^u$ is unitarily equivalent to $N^v$.  Thus
we know that $[u]=[v]$.
\end{proof}

\section{Main Theorem}

Now we have the continuous injection we need to prove the main
result.

\begin{theorem}
\label{genccr}
Suppose $G$ is groupoid in which all of the
isotropy groups are amenable.
Then $\cs(G)$ is CCR if and only
if $\go/G$ is $T_1$ and each of the isotropy groups are CCR.
\end{theorem}

\begin{proof}

Suppose $\cs(G)$ is CCR.  This means that points of
$\cs(G)^\wedge$ are closed.  Because we have a continuous
injection between $\go/G$ and $\cs(G)^\wedge$,
we know points of $\go/G$ are also
closed.  Thus $\go/G$ is $T_1$.

We know that every representation of
$\cs(G)$
factors through $\csg{\overline{[u]}} = \csg{[u]}$ for some
$u\in \go$ from \cite[Corollary~3.5]{mypaper}.
But $\csg{[u]}$ is a transitive groupoid thus
$\csg{[u]} \cong \cs(G^u_u) \otimes K$
\cite[Theorem~3.1]{equivalence}.  This means that any representation of $\cs(G^u_u)$ lifts to a representation of a CCR algebra, namely  $\cs(G)$.  Since the lifted representation is onto the compact operators, so must the original have been.  That is, each of the isotropy groups are CCR.

Now, conversely, suppose $\go/G$ is $T_1$ and isotropy groups are CCR.  Again,
we know from \cite[Corollary~3.5]{mypaper} that every representation $L$ of
$\cs(G)$
factors through $\csg{\overline{[u]}} = \csg{[u]}$ for some
$u\in \go$ and that $\csg{[u]} \cong \cs(G^u_u) \otimes K$ by
\cite[Theorem~3.1]{equivalence}.  This is CCR
because
$\cs(G^u_u)$ is assumed to be CCR.  This means that $L$  is
lifted from a representation  of a CCR $\cs$-algebra making
$L$ a representation onto the compact operators.  That is, $\cs(G)$ is CCR.

\end{proof}

\section{GCR Result}
\begin{theorem}
\label{genGCR}
Suppose $G$ is groupoid in which all of the isotropy groups are amenable, then $\cs(G)$ is
GCR if and only
if $\go/G$ is $T_0$ and all of the isotropy groups are GCR.
\end{theorem}

\begin{proof}
Suppose $\cs(G)$ is GCR. This means that $\cs(G)^\wedge$ is $T_0$.
Since there is a continuous injection
from the orbit space to the spectrum, the orbit space must also
be $T_0$.

From \cite[Proposition~5.1]{mypaper}, we
know that every irreducible representation of
$\cs(G)$ is the canonical extension of a representation of
$\cs(G_{U_\alpha \setminus U_{\alpha-1}})$ where
$U_\alpha \setminus U_{\alpha-1}$ is Hausdorff.  Following the proof in the CCR case,  we know that every representation of
$\cs(G_{U_\alpha \setminus U_{\alpha-1}})$
factors through $\csg{\overline{[u]}} = \csg{[u]}$ for some
$u\in U_\alpha \setminus U_{\alpha-1}$ and that $\csg{[u]} \cong \cs(G^u_u) \otimes K$ by
\cite[Theorem~3.1]{equivalence}.  Now if we consider a representation $l$ of $\cs(G^u_u)$ for any $u \in \go$, it lifts to a representation of $\cs(G_{U_\alpha \setminus U_{\alpha-1}})$ which can be extened canonically to $\cs(G)$.  Since the later is GCR, then so is $l$.

Conversely, suppose that $\go/G$ is $T_0$ and the isotropy groups are GCR.   Once again, we
know that every irreducible representation $L$ of
$\cs(G)$ is the canonical extension of a representation of
$\cs(G_{U_\alpha \setminus U_{\alpha-1}})$ where
$U_\alpha \setminus U_{\alpha-1}$ is Hausdorff.  Also,
we know that every representation of
$\cs(G_{U_\alpha \setminus U_{\alpha-1}})$
factors through $\csg{\overline{[u]}} = \csg{[u]}$ for some
$u\in U_\alpha \setminus U_{\alpha-1}$ and that $\csg{[u]} \cong \cs(G^u_u) \otimes K$ by
\cite[Theorem~3.1]{equivalence}.  Since $\cs(G^u_u)$ is assumed to be GCR, then $\cs(G)$ must also be GCR.

\end{proof}

\thanks{This research was done as part of the
author's Ph.D. thesis under the direction of Dana~P.
Williams.  Thank you to Dana for his continued support. }

\end{document}